\documentclass[12pt]{article}
\usepackage{amscd}
\usepackage{amsmath,amsfonts,amssymb,amscd}
\usepackage{indentfirst,graphics,epsfig,psfrag,cite}
\usepackage{ifpdf}
\input{epsf}

\makeatletter \@addtoreset{figure}{section} \makeatother
\makeatletter
\long\def\@makecaption#1#2{%
   \vskip 10\p@
   \setbox\@tempboxa\hbox{{#1}\ \ #2}%
   \ifdim \wd\@tempboxa >\hsize
       {#1}\ \ #2\par
   \else
       \hbox to\hsize{\hfil\box\@tempboxa\hfil}%
   \fi}
\makeatother

\newtheorem{theo}{Theorem}[section]
\newtheorem{cor}[theo]{Corollary}
\newtheorem{lem}[theo]{Lemma}

\newcommand{\qed}{{\hfill\rule{3pt}{7pt}}}
\def\pf{\noindent {\it Proof.} }

\def\qed{\hfill \rule{4pt}{7pt}}
\def\pf{\noindent {\it Proof.} }

\begin{document}

\begin{center} {\Large \bf On the maximal energy tree\\[2mm]
with two maximum degree vertices}\footnote{Supported by NSFC and
``the Fundamental Research Funds for the Central Universities". }

\end{center}
\pagestyle{plain}
\begin{center}
{
  {\small  Jing Li,  Xueliang Li, Yongtang Shi}\\[3mm]
  {\small Center for Combinatorics and LPMC-TJKLC}\\
  {\small Nankai University, Tianjin 300071, China}\\
  {\small E-mail: lj02013@163.com; lxl@nankai.edu.cn; shi@nankai.edu.cn}\\[2mm]

}
\end{center}

\begin{center}
\begin{minipage}{140mm}
\vskip 0.3cm
\begin{center}

{\bf Abstract}
\end{center}

{\small  For a simple graph $G$, the energy $E(G)$ is defined as the
sum of the absolute values of all eigenvalues of its adjacent
matrix. For $\Delta\geq 3$ and $t\geq 3$, denote by $T_a(\Delta,t)$
(or simply $T_a$) the tree formed from a path $P_t$ on $t$ vertices
by attaching $\Delta-1$ $P_2$'s on each end of the path $P_t$, and
$T_b(\Delta, t)$ (or simply $T_b$) the tree formed from $P_{t+2}$ by
attaching $\Delta-1$ $P_2$'s on an end of the $P_{t+2}$ and $\Delta
-2$ $P_2$'s on the vertex next to the end. In [X. Li, X. Yao, J.
Zhang and I. Gutman, Maximum energy trees with two maximum degree
vertices, J. Math. Chem. 45(2009), 962--973], Li et al. proved that
among trees of order $n$ with two vertices of maximum degree
$\Delta$, the maximal energy tree is either the graph $T_a$ or the
graph $T_b$, where $t=n+4-4\Delta\geq 3$. However, they could not
determine which one of $T_a$ and $T_b$ is the maximal energy tree.
This is because the quasi-order method is invalid for comparing
their energies. In this paper, we use a new method to determine the
maximal energy tree. It turns out that things are more complicated.
We prove that the maximal energy tree is $T_b$ for $\Delta\geq 7$
and any $t\geq 3$, while the maximal energy tree is $T_a$ for
$\Delta=3$ and any $t\geq 3$. Moreover, for $\Delta=4$, the maximal
energy tree is $T_a$ for all $t\geq 3$ but $t=4$, for which $T_b$ is
the maximal energy tree. For $\Delta=5$, the maximal energy tree is
$T_b$ for all $t\geq 3$ but $t$ is odd and $3\leq t\leq 89$, for
which $T_a$ is the maximal energy tree. For $\Delta=6$, the maximal
energy tree is $T_b$ for all $t\geq 3$ but $t=3,5,7$, for which
$T_a$ is the maximal energy tree. One can see that for most
$\Delta$, $T_b$ is the maximal energy tree, $\Delta=5$ is a turning
point, and $\Delta=3$ and 4 are exceptional cases.}\\[2mm]
{\bf Keywords:} graph energy, tree, Coulson integral formula.\\[2mm]
{\bf AMS subject classification 2010:} 05C50, 05C90, 15A18, 92E10.

\end{minipage}
\end{center}

\section{Introduction}

Let $G$ be a simple graph of order $n$, it is well known \cite{Cvet}
that the characteristic polynomial of $G$ has the form
$$\varphi(G,x)=\sum_{k=0}^na_kx^{n-k}.$$
The match polynomial of $G$ is defined as
$$m(G,x)=\sum_{k=0}^{\lfloor n/2 \rfloor} (-1)^km(G,k)x^{n-2k},$$
where $m(G,k)$ denotes the number of $k$-matchings of $G$ and
$m(G,0)=1$. If $G=T$ is a tree of order $n$, then
$$\varphi(T,x)=m(T,x)=\sum_{k=0}^{\lfloor n/2 \rfloor} (-1)^km(T,k)x^{n-2k}.$$

Let $\lambda_1,\lambda_2,\cdots,\lambda_n$ be the eigenvalues of
$G$, then the energy of $G$ is defined as
$$E(G)=\sum_{i=1}^n|\lambda_i|,$$ which was introduced by Gutman in \cite{Gut1}.
If $T$ is a tree of order $n$, then by Coulson integral formula
\cite{Gut2, Gut3}, we have
$$E(T)=\frac{2}{\pi}\int_0^{+\infty}\frac{1}{x^2}\log \left[ \sum_{k=0}
^{\lfloor n/2 \rfloor}m(T,k)x^{2k} \right]dx.$$ In order to avoid
the signs in the matching polynomial, this immediately motivates us
to introduce a new graph polynomial
$$m^+(G,x)=\sum_{k=0}^{\lfloor n/2 \rfloor}m(G,k)x^{2k}.$$
Then we have
\begin{equation}\label{LJequ1}
E(T)=\frac{2}{\pi}\int_0^{+\infty}\frac{1}{x^2}\log m^+(T,x)dx.
\end{equation}
Although $m^+(G,x)$ is nothing new but
$m^+(G,x)=(ix)^nm(G,(ix)^{-1})$, we shall see later that this will
bring us a lot of computational convenience. Some basic properties
of $m^+(G,x)$ will be given in next section.

We refer to the survey \cite{gutman&lxl2009} for more results on
graph energy. For terminology and notation not defined here, we
refer to the book of Bondy and Murty \cite{Bondy}.

\begin{figure}[ht]
\begin{center}
\includegraphics[width=15cm]{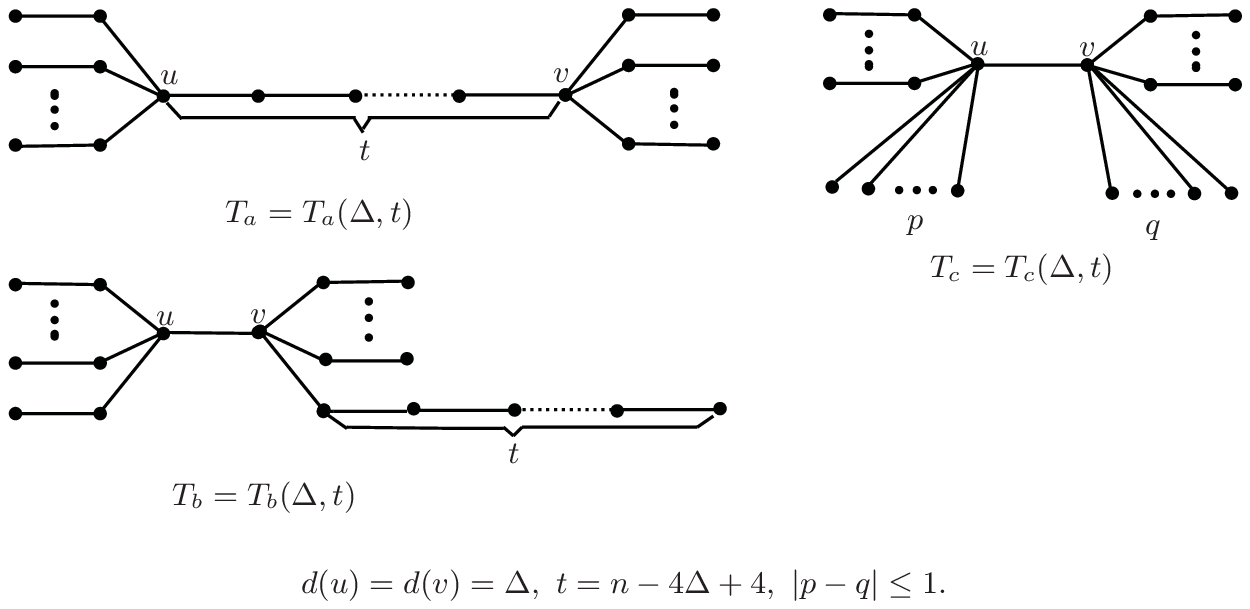}
\end{center}
\caption{The maximal energy trees with $n$ vertices and two vertices
$u$, $v$ of maximum degree $\Delta$.}\label{LJfig1}
\end{figure}

Graphs with extremal energies are interested in literature. Gutman
\cite{Gut2} proved that the star and the path has the minimal and
the maximal energy among all trees, respectively. Lin et al.
\cite{WLin} showed that among trees with a fixed number of vertices
($n$) and of maximum vertex degree ($\Delta$), the maximal energy
tree has exactly one branching vertex (of degree $\Delta$) and as
many as possible 2-branches. Li et al. \cite{XLi} gave the following
Theorem \ref{thmli} about the maximal energy tree with two maximum
degree vertices. In a similar way, Yao \cite{Yao} studied the
maximal energy tree with one maximum and one second maximum degree
vertex. A {\it branching vertex} is a vertex whose degree is three
or greater, and a pendent vertex attached to a vertex of degree two
is called a {\it $2$-branch}.
\begin{theo}[\cite{XLi}] \label{thmli} Among trees with a fixed number of
vertices $(n)$ and two vertices of maximum  degree $(\Delta)$, the
maximal energy
tree has as many as possible $2$-branches.\\
(1) If $n\leq 4\Delta-2$, then the maximal energy tree is the graph
$T_c=T_c(\Delta, t)$ depicted in Figure \ref{LJfig1}, in which the
numbers of pendent vertices attached to the two branching vertices
$u$ and $v$
differ by at most 1. \\
(2) If $n\geq 4\Delta-1$, then the maximal energy tree is either the
graph $T_a=T_a(\Delta, t)$ or the graph $T_b=T_b(\Delta, t)$,
depicted in Figure \ref{LJfig1}. \label{theo1}
\end{theo}

From Theorem \ref{theo1}, one can see that for $n\geq 4\Delta-1$,
they could not determine which one of the graphs $T_a$ and $T_b$ has
the maximal energy. They gave small examples showing that both cases
could happen. In fact, the quasi-order method they used before is
invalid for the special case. Recently, for these quasi-order
incomparable problems, Huo et al. found an efficient way to
determine which one attains the extremal value of the energy, we
refer to \cite{Huo1, Huo2, Huo3, Huo7, Huo4, Huo5, Huo6} for
details. In this paper, we will use this newly developed method to
determine which one of the graphs $T_a$ and $T_b$ has the maximal
energy, solving this unsolved problem. It turns that this problem is
more complicated than those in \cite{Huo1, Huo2, Huo3, Huo7, Huo4,
Huo5, Huo6}.

\section{Preliminaries}

In this section, we will give some properties of the new polynomial
$m^+(G,x)$, which will be used in what follows. The proofs are
omitted, since they are the same as those for matching polynomial.

\begin{lem}\label{LJlem1}
Let $K_n$ be a complete graph with $n$ vertices and $\overline{K_n}$
the complement of $K_n$, then
$$m^+(\overline{K_n},x)=1,$$
for any $n\geq 0$, defining $m^+(\overline{K_0},x)=1$, where both
$K_0$ and $\overline{K_0}$ are the null graph.
\end{lem}

Similar to the properties of matching polynomial, we have
\begin{lem}\label{LJlem2}
Let $G_1$ and $G_2$ be two vertex disjoint graphs. Then
$$m^+(G_1\cup G_2,x)=m^+(G_1,x)\cdot m^+(G_2,x).$$
\end{lem}

\begin{lem}\label{LJlem3}
Let $e=uv$ be an edge of graph $G$. Then we have
$$m^+(G,x)=m^+(G-e,x)+x^2 m^+(G-u-v,x).$$
\end{lem}

\begin{lem}\label{LJlem4}
Let $v$ be a vertex of $G$ and $N(v)=\{v_1,v_2,\ldots,v_r\}$ the set
of all neighbors of $v$ in $G$. Then
\begin{equation*}
m^+(G,x)=m^+(G-v,x)+x^2 \sum_{v_i\in N(v)}m^+(G-v-v_i,x).
\end{equation*}
\end{lem}

The following recursive equations can be gotten from Lemma
\ref{LJlem3} immediately.
\begin{lem}\label{LJlem5}
Let $P_t$ denote a path on $t$ vertices. Then
\begin{eqnarray*}
(1)\ m^+(P_t,x)&=&m^+(P_{t-1},x)+x^2
m^+(P_{t-2},x),\ \mbox{for any~} t\geq 1,\\
(2)\ m^+(P_t,x)&=&(1+x^2)m^+(P_{t-2},x)+x^2 m^+(P_{t-3},x),\
\mbox{for any~} t\geq 2.
\end{eqnarray*}
The initials are $m^+(P_0,x)=m^+(P_1,x)=1$, and we define
$m^+(P_{-1},x)=0$.
\end{lem}

From Lemma \ref{LJlem5}, one can easily obtain

\begin{cor}\label{LJcor1}
Let $P_t$ be a path on $t$ vertices. Then for any real number $x$,
\begin{equation*}
m^+(P_{t-1},x)\leq m^+(P_{t},x)\leq (1+x^2)m^+(P_{t-1},x),
 \ \mbox{for any~} t \geq 1.
 \end{equation*}
\end{cor}
Although $m^+(G,x)$ has many other properties, the above ones are
enough for our use.

\section{Main results}

Before giving our main results, we state some knowledge on real
analysis, for which we refer to \cite{Zori}.

\begin{lem}\label{LJlem6}
For any real number $X>-1$, we have
\begin{equation*}
\frac{X}{1+X}\leq \log(1+X)\leq X.
\end{equation*}
\end{lem}

To compare the energies of $T_a$ and $T_b$, or more precisely,
$T_a(\Delta, t)$ and $T_b(\Delta, t)$, means to compare the values
of two functions with the parameters $\Delta$ and $t$, which are
denoted by $E(T_a(\Delta,t))$ and $E(T_b(\Delta,t))$. Since
$E(T_a(2,t))=E(T_b(2,t))$ for any $t\geq 2$ and
$E(T_a(\Delta,2))=E(T_b(\Delta,2))$ for any $\Delta\geq 2$, we
always assume that $\Delta\geq 3$ and $t\geq 3$.

For notational convenience, we introduce the following things:
$$A_1=(1+x^2)(1+\Delta x^2)(2x^4+(\Delta+2)x^2+1),$$
$$A_2=x^2(1+x^2)(x^6+(\Delta^2+2)x^4+(2\Delta+1)x^2+1),$$
$$B_1=(\Delta+2)x^8+(2\Delta^2+6)x^6+(\Delta^2+4\Delta+4)
x^4+(2\Delta+3)x^2+1,$$
$$B_2=x^2(1+x^2)(x^6+(\Delta^2+2)x^4+(2\Delta+1)x^2+1).$$ Using Lemmas
\ref{LJlem4} and \ref{LJlem5} repeatedly, we can easily get the
following two recursive formulas:
\begin{eqnarray}\label{LJequ8}
m^+(T_a,x)=(1+x^2)^{2\Delta-5}(A_1m^+(P_{t-3},x)+A_2m^+(P_{t-4},x)),
\end{eqnarray}
and
\begin{eqnarray}\label{LJequ9}
m^+(T_b,x)=(1+x^2)^{2\Delta-5}(B_1m^+(P_{t-3},x)+B_2m^+(P_{t-4},x)),
\end{eqnarray}
From Eqs. \eqref{LJequ8} and \eqref{LJequ9}, by some elementary
calculations we can obtain
\begin{eqnarray}\label{LJequ10}
m^+(T_a,x)-m^+(T_b,x)=(1+x^2)^{2\Delta-5}(\Delta-2)x^6
(x^2-(\Delta-2))m^+(P_{t-3},x).
\end{eqnarray}

Now we give one of our main results.
\begin{theo}\label{LJtheo1}
Among trees with $n$ vertices and two vertices of maximum degree
$\Delta$, the maximal energy tree has as many as possible
2-branches. If $\Delta\geq 8$ and $t\geq 3$, then the maximal energy
tree is the graph $T_b$, where $t=n+4-4\Delta$.
\end{theo}
\pf From Eq. (\ref{LJequ1}), we have
\begin{eqnarray}
E(T_a)-E(T_b)&=&\frac{2}{\pi}\int_0^{+\infty}\frac{1}{x^2}\log
\frac{m^+(T_a,x)}{m^+(T_b,x)}dx\nonumber\\[2mm]
&=&\frac{2}{\pi}\int_0^{+\infty}\frac{1}{x^2}\log
\left(1+\frac{m^+(T_a,x)-m^+(T_b,x)}{m^+(T_b,x)}\right)dx.\label{LJequ11}
\end{eqnarray}
We express $g(\Delta,t,x)$ as
$$g(\Delta,t,x)=\frac{1}{x^2}\log\left(1+\frac{m^+(T_a,x)-m^+(T_b,x)} {m^+(T_b,x)}\right).$$
Since $m^+(T_a,x)> 0$ and $m^+(T_b,x)>0$, we have
\begin{equation*}
\frac{m^+(T_a,x)-m^+(T_b,x)}{m^+(T_b,x)}=\frac{m^+(T_a,x)}{m^+(T_b,x)}-1>
-1.
\end{equation*}
Therefore, by Lemma \ref{LJlem6} we have
\begin{equation}\label{LJequ7}
\frac{1}{x^2}\cdot\frac{m^+(T_a,x)-m^+(T_b,x)}{m^+(T_a,x)}\leq
g(\Delta,t,x)\leq
\frac{1}{x^2}\cdot\frac{m^+(T_a,x)-m^+(T_b,x)}{m^+(T_b,x)}.
\end{equation}
Substituting the recursive formulas (\ref{LJequ8}), (\ref{LJequ9})
and (\ref{LJequ10}) to Eq. (\ref{LJequ7}), we get that
\begin{eqnarray*}
g(\Delta,t,x)
 &\leq&\frac{1}{x^2}\cdot\frac{(1+x^2)^{2\Delta-5}(\Delta-2)x^6
(x^2-(\Delta-2))m^+(P_{t-3},x)}{(1+x^2)^{2\Delta-5}
(B_1m^+(P_{t-3},x)+B_2m^+(P_{t-4},x))} \\[2mm]
&=&\frac{(\Delta-2)x^4
(x^2-(\Delta-2))m^+(P_{t-3},x)}{B_1m^+(P_{t-3},x)+B_2m^+(P_{t-4},x)},
\end{eqnarray*}
and
\begin{eqnarray*}
g(\Delta,t,x)
 &\geq&\frac{1}{x^2}\cdot\frac{(1+x^2)^{2\Delta-5}(\Delta-2)x^6
(x^2-(\Delta-2))m^+(P_{t-3},x)}{(1+x^2)^{2\Delta-5}
(A_1m^+(P_{t-3},x)+A_2m^+(P_{t-4},x))} \\[2mm]
&=&\frac{(\Delta-2)x^4
(x^2-(\Delta-2))m^+(P_{t-3},x)}{A_1m^+(P_{t-3},x)+A_2m^+(P_{t-4},x)}.
\end{eqnarray*}
By Corollary \ref{LJcor1}, we have $m^+(P_{t-4},x)\leq
m^+(P_{t-3},x)$ and $m^+(P_{t-4},x)\geq
\frac{m^+(P_{t-3},x)}{1+x^2}$ for $\Delta\geq 3$ and $t\geq 4$. Then
if $x\geq \sqrt{\Delta-2}$,
\begin{eqnarray*}
|g(\Delta,t,x)|
 &\leq&\frac{(\Delta-2)x^4
(x^2-(\Delta-2))}{
B_1+B_2/(1+x^2)}\\[2mm]
&=&\frac{(\Delta-2)x^4
(x^2-(\Delta-2))}{(\Delta+3)x^8+(3\Delta^2+8)x^6+(\Delta^2+6\Delta+5)x^4
+(2\Delta+4)x^2+1},
\end{eqnarray*}
and if $x\leq \sqrt{\Delta-2}$,
\begin{eqnarray*}
|g(\Delta,t,x)|
 &\leq&\frac{(\Delta-2)x^4
(\Delta-2-x^2)}{
A_1+A_2/(1+x^2)}\\[2mm]
&=&\frac{(\Delta-2)x^4
(\Delta-2-x^2)}{(2\Delta+1)x^8+(2\Delta^2+4\Delta+4)x^6+
(\Delta^2+6\Delta+5)x^4+(2\Delta+4)x^2+1}.
\end{eqnarray*}
Since for $\Delta\geq 3$ and any $x\ge 0$, we always have
\begin{equation*}
(\Delta-2)x^4 (x^2-(\Delta-2))(1+x^2)\leq
(\Delta+3)x^8+(3\Delta^2+8)x^6+(\Delta^2+6\Delta+5)x^4
+(2\Delta+4)x^2+1,
\end{equation*}
and
\begin{equation*}
(\Delta-2)x^4 (\Delta-2-x^2)(1+x^2)\leq
(2\Delta+1)x^8+(2\Delta^2+4\Delta+4)x^6+
(\Delta^2+6\Delta+5)x^4+(2\Delta+4)x^2+1,
\end{equation*}
we can get that for $\Delta\geq 3$ and any $x\ge 0$,
$$|g(\Delta,t,x)|
 \leq\frac{1}{1+x^2},$$
 while $\int_0^{+\infty}\frac{2}{1+x^2}dx=\frac{\pi}{2}$ is
 convergent. From the well-known Weierstrass's criterion
 (for example, see \cite{Zori}), we can get that $E(T_a)-E(T_b)
 =\frac{2}{\pi}\int_0^{+\infty}g(\Delta,t,x)dx$ is uniformly
 convergent. Then
 \begin{equation*}
\frac{2}{\pi}\int_0^{+\infty}\frac{1}{x^2}\cdot\frac{m^+(T_a,x)-m^+(T_b,x)}
{m^+(T_a,x)}dx\leq E(T_a)-E(T_b)\leq
\frac{2}{\pi}\int_0^{+\infty}\frac{1}{x^2}\cdot\frac{m^+(T_a,x)-m^+(T_b,x)}
{m^+(T_b,x)}dx.
\end{equation*}
Thus, for $t\geq 4$, we have
\begin{eqnarray*}
&&E(T_a)-E(T_b)\\[2mm]
&\leq&\frac{2}{\pi}\int_0^{+\infty}\frac{1}{x^2}\cdot
\frac{m^+(T_a,x)-m^+(T_b,x)}{m^+(T_b,x)}dx\\[2mm]
&=&\frac{2}{\pi}\int_0^{+\infty}\frac{(\Delta-2)x^4
(x^2-(\Delta-2))m^+(P_{t-3},x)}{
B_1m^+(P_{t-3},x)+B_2m^+(P_{t-4},x)}dx\\[2mm]
&\leq&\frac{2}{\pi}\int_{\sqrt{\Delta-2}}^{+\infty}
\frac{(\Delta-2)x^4 (x^2-(\Delta-2))}{B_1+\frac{B_2}{1+x^2}}dx
-\frac{2}{\pi}\int_0^{\sqrt{\Delta-2}} \frac{(\Delta-2)x^4
(\Delta-2-x^2)}{B_1+B_2}dx.
\end{eqnarray*}
We calculate the two parts respectively. The first part is
\begin{eqnarray*}
&&\frac{2}{\pi}\int_{\sqrt{\Delta-2}}^{+\infty} \frac{(\Delta-2)x^4
(x^2-(\Delta-2))}{B_1+\frac{B_2}{1+x^2}}dx\\[2mm]
&=&\frac{2}{\pi}\int_{\sqrt{\Delta-2}}^{+\infty} \frac{(\Delta-2)x^4
(x^2-(\Delta-2))}{(\Delta+3)x^8+(3\Delta^2+8)x^6+(\Delta^2+6\Delta+5)x^4
+(2\Delta+4)x^2+1}dx\\[2mm]
&<&\frac{2}{\pi}\int_{\sqrt{\Delta-2}}^{+\infty} \frac{(\Delta-2)x^4
(x^2-(\Delta-2))}{(\Delta+3)x^8}dx=\frac{2}{\pi}\cdot
\frac{2\sqrt{\Delta-2}}{3(\Delta+3)}.
\end{eqnarray*}
The second part is
\begin{eqnarray*}
&&\frac{2}{\pi}\int_0^{\sqrt{\Delta-2}} \frac{(\Delta-2)x^4
(\Delta-2-x^2)}{B_1+B_2}dx\\[2mm]
&=&\frac{2}{\pi}\int_0^{\sqrt{\Delta-2}} \frac{(\Delta-2)x^4
(\Delta-2-x^2)}{h(\Delta, x)}dx\\[2mm]
&>&\frac{2}{\pi}\int_0^1 \frac{(\Delta-2)x^4
(\Delta-2-x^2)}{\frac{5\Delta^2+11\Delta+26}{2}(x^2+1)}dx+\frac{2}{\pi}\int_1
^{\sqrt{\Delta-2}} \frac{(\Delta-2)x^4
(\Delta-2-x^2)}{(5\Delta^2+11\Delta+26)x^{10}}dx\\[2mm]
&=&\frac{2}{\pi}\left(\frac{-45\pi\Delta-34\Delta^2+74\Delta+30\pi-12+15\pi\Delta
^2+\frac{4}{\sqrt{\Delta-2}}}{30(26+11\Delta+5\Delta^2)}\right),
\end{eqnarray*}
where $h(\Delta,
x)=x^{10}+(\Delta^2+\Delta+5)x^8+(3\Delta^2+2\Delta+9)x^6+
(\Delta^2+6\Delta+6)x^4+(2\Delta+4)x^2+1$. Now, when $\Delta\geq
65$, we have that
\begin{eqnarray*}
&&E(T_a)-E(T_b)\\[2mm]
&<&\frac{2}{\pi}\cdot \frac{2\sqrt{\Delta-2}}{3(\Delta+3)}-
\frac{2}{\pi}\left(\frac{-45\pi\Delta-34\Delta^2+74\Delta+30\pi-12+15\pi\Delta
^2+\frac{4}{\sqrt{\Delta-2}}}{30(26+11\Delta+5\Delta^2)}\right)\leq0.
\end{eqnarray*}
For $t=3$, we have $m^+(P_{t-4},x)=m^+(P_{-1},x)=0$. By a similar
method as above, we can get that $E(T_a)-E(T_b)<0$ when $\Delta\geq
24$.

Therefore, for $\Delta\geq 65$ and $t\geq3$, we have
$E(T_a)<E(T_b)$.

For $8\leq \Delta\leq 64$, we can calculate
\begin{equation*}
E(T_a)-E(T_b)\leq \frac{2}{\pi}\cdot f(\Delta,x)<0
\end{equation*}
directly by computer programm, as shown in Table \ref{LJtab1}, where
\begin{equation*}
f(\Delta,x)=\int_{\sqrt{\Delta-2}}^{+\infty} \frac{(\Delta-2)x^4
(x^2-(\Delta-2))}{B_1+\frac{B_2}{1+x^2}}dx -\int_0^{\sqrt{\Delta-2}}
\frac{(\Delta-2)x^4 (\Delta-2-x^2)}{B_1+B_2}dx.
\end{equation*}

\begin{table}[ht]
\centering
\begin{tabular}{|cc|cc|cc|cc|}
  \hline
  $\Delta$& $f(\Delta,x)$ &$\Delta$ & $f(\Delta,x)$ &$\Delta$ &
  $f(\Delta,x)$ &$\Delta$ & $f(\Delta,x)$\\
  \hline
  8 & -0.00377 &  23 & -0.20792 & 38 & -0.29961 &  53 & -0.35353 \\
  9 & -0.02418 &  24 & -0.21611 & 39 & -0.30403 &  54 & -0.35638 \\
  10 & -0.04352 & 25 & -0.22390 & 40 & -0.30830 &  55 & -0.35917 \\
  11 & -0.06168 & 26 & -0.23132 & 41 & -0.31244 &  56 & -0.36188 \\
  12 &-0.07866 &  27 & -0.23841 & 42 & -0.31644 &  57 & -0.36454 \\
  13 & -0.09452 & 28 & -0.24518 & 43 & -0.32032 &  58 & -0.36713 \\
  14 & -0.10933 & 29 & -0.25165 & 44 & -0.32409 &  59 & -0.36965 \\
  15 & -0.12317 & 30 & -0.25786 & 45 & -0.32774 &  60 & -0.37213 \\
  16 & -0.13613 & 31 & -0.26381 & 46 & -0.33129 &  61 & -0.37454 \\
  17 & -0.14829 & 32 & -0.26953 & 47 & -0.33473 &  62 & -0.37691 \\
  18 & -0.15972 & 33 & -0.27502 & 48 & -0.33808 &  63 & -0.37922 \\
  19 & -0.17048 & 34 & -0.28031 & 49 & -0.34134 &  64 & -0.38148 \\
  20 & -0.18063 & 35 & -0.28540 & 50 & -0.34451 &  65 & -0.38369 \\
  21 & -0.19022 & 36 & -0.29031 & 51 & -0.34759 &  66 & -0.38586 \\
  22 & -0.19931 & 37 & -0.29504 & 52 & -0.35060 &  67 & -0.38798 \\
  \hline
\end{tabular}
  \caption{The values of $f(\Delta,x)$ for $8\leq \Delta\leq 67$.}\label{LJtab1}
\end{table}

The proof is thus complete.\qed

Now we are left with the cases $3\leq \Delta \leq 7$. At first, we
consider the case of $\Delta=3$ and $t\geq 3$. In this case, we have
$n=4\Delta-4+t\geq 11$.

\begin{theo}\label{LJtheo2}
Among trees with $n$ vertices and two vertices of maximum degree
$\Delta=3$, the maximal energy tree has as many as possible
$2$-branches. If $n\geq 11$, then the maximal energy tree is the
graph $T_a$.
\end{theo}
\pf For $\Delta=3$ and $t\geq 4$, by Eqs. (\ref{LJequ1}),
(\ref{LJequ7}) and Corollary\ref{LJcor1}, we have
\begin{eqnarray*}
E(T_a)-E(T_b)&\geq&\frac{2}{\pi}\int_0^{+\infty}\frac{1}{x^2}\cdot
\frac{m^+(T_a,x)-m^+(T_b,x)}{m^+(T_a,x)}dx\\[2mm]
&=&\frac{2}{\pi}\int_0^{+\infty}\frac{1}{x^2}\cdot \frac{x^6
(x^2-1)m^+(P_{t-3},x)}{A_1m^+(P_{t-3},x)+A_2m^+(P_{t-4},x)}dx\\[2mm]
&\geq&\frac{2}{\pi}\int_1^{+\infty} \frac{x^4 (x^2-1)}{A_1+A_2}dx
-\frac{2}{\pi}\int_0^1 \frac{x^4
(1-x^2)}{A_1+\frac{A_2}{1+x^2}}dx\\[2mm]
&=&\frac{2}{\pi}\int_1^{+\infty} \frac{x^4
(x^2-1)}{x^{10}+18x^8+41x^6+33x^4+10x^2+1}dx\\[2mm]
&& -\frac{2}{\pi}\int_0^1 \frac{x^4
(1-x^2)}{7x^8+34x^6+32x^4+10x^2+1}dx\\[2mm]
&>&\frac{2}{\pi}\cdot 0.00996>0.
\end{eqnarray*}

For $\Delta=3$ and $t=3$, we can compare the energies of the two
graphs directly and get that $E(T_a)>E(T_b)$.

Therefore, for $\Delta=3$ and $t\geq 3$, we have
$E(T_a)>E(T_b)$.\qed

Now we give two lemmas about the properties of the new polynomial
$m^+(P_t,x)$.

\begin{lem}\label{LJlem3.4}
For $t\geq -1$, the polynomial $m^+(P_t,x)$ has the following form
\begin{equation*}
m^+(P_t,x)=\frac{1}{\sqrt{1+4x^2}}
(\lambda_1^{t+1}-\lambda_2^{t+1}),
\end{equation*}
where $\lambda_1=\frac{1+\sqrt{1+4x^2}}{2}$ and
$\lambda_2=\frac{1-\sqrt{1+4x^2}}{2}$.
\end{lem}
\pf By Lemma \ref{LJlem5}, $m^+(P_t,x)=m^+(P_{t-1},x)+x^2
m^+(P_{t-2},x)\ \mbox{for any~} t\geq 1$. Thus, it satisfies the
recursive formula $h(t,x)=h(t-1,x)+x^2h(t-2,x)$, and the  general
solution of this linear homogeneous recurrence relation is
$h(t,x)=P(x)\lambda_1^t+Q(x)\lambda_2^t$, where
$\lambda_1=\frac{1+\sqrt{1+4x^2}}{2}$ and
$\lambda_2=\frac{1-\sqrt{1+4x^2}}{2}$. Considering the initial
values $m^+(P_1,x)=1$ and $m^+(P_2,x)=1+x^2$, by some elementary
calculations, we can easily obtain that
\begin{eqnarray*}
&P(x)=\frac{1+\sqrt{1+4x^2}}{2\sqrt{1+4x^2}},\qquad
Q(x)=\frac{-1+\sqrt{1+4x^2}}{2\sqrt{1+4x^2}}.&
\end{eqnarray*}
Thus,
\begin{eqnarray*}
m^+(P_t,x)=P(x)\lambda_1^t+Q(x)\lambda_2^t=\frac{1}{\sqrt{1+4x^2}}
(\lambda_1^{t+1}-\lambda_2^{t+1}).
\end{eqnarray*}

As we have defined, the initials are $m^+(P_{-1},x)=0$ and
$m^+(P_0,x)=1$, from which we can get the result for all $t\geq
-1$.\qed

\begin{lem}\label{LJlem3.5}
Suppose $t\geq 4$. If  $t$ is even, then
\begin{equation*}
\frac{2}{1+\sqrt{1+4x^2}}<\frac{m^+(P_{t-4},x)}{m^+(P_{t-3},x)}\leq
1.
\end{equation*}
 If $t$ is odd, then
\begin{equation*}
\frac{1}{1+x^2}\leq\frac{m^+(P_{t-4},x)}{m^+(P_{t-3},x)}<
\frac{2}{1+\sqrt{1+4x^2}}.
\end{equation*}
\end{lem}
\pf From Corollary \ref{LJcor1}, we know that
\begin{equation*}
\frac{1}{1+x^2}\leq\frac{m^+(P_{t-4},x)}{m^+(P_{t-3},x)}\leq 1.
\end{equation*}
By the definitions of $\lambda_1$ and $\lambda_2$, we conclude that
$\lambda_1>0$ and $\lambda_2<0$ for any $x$. By Lemma
\ref{LJlem3.4}, if $t$ is even, then
\begin{eqnarray*}
\frac{m^+(P_{t-4},x)}{m^+(P_{t-3},x)}- \frac{2}{1+\sqrt{1+4x^2}}
=\frac{\lambda_1^{t-3}-\lambda_2^{t-3}}{\lambda_1^{t-2}-
\lambda_2^{t-2}}-\frac{1}{\lambda_1}
=\frac{-\lambda_2^{t-3}(\lambda_1-\lambda_2)}
{\lambda_1(\lambda_1^{t-2}- \lambda_2^{t-2})}>0.
\end{eqnarray*}
Thus,
\begin{equation*}
\frac{2}{1+\sqrt{1+4x^2}}<\frac{m^+(P_{t-4},x)}{m^+(P_{t-3},x)}\leq
1.
\end{equation*}
If $t$ is odd, then obviously
\begin{equation*}
\frac{1}{1+x^2}\leq\frac{m^+(P_{t-4},x)}{m^+(P_{t-3},x)}<
\frac{2}{1+\sqrt{1+4x^2}}.
\end{equation*} \qed

Now we deal with the case $\Delta=4$ and $t\geq 3$.

\begin{theo}\label{LJtheo3}
Among trees with $n$ vertices and two vertices of maximum degree
$\Delta=4$, the maximal energy tree has as many as possible
$2$-branches. The maximal energy tree is the graph $T_b$ if $t=4$,
and the graph $T_a$ otherwise, where $t=n+4-4\Delta$.
\end{theo}
\pf  By Eqs. (\ref{LJequ8}), (\ref{LJequ9}), (\ref{LJequ10}) and
(\ref{LJequ11}), we have
\begin{eqnarray}\label{LJequ2011}
E(T_a)-E(T_b)&=&\frac{2}{\pi}\int_0^{+\infty}\frac{1}{x^2}\log\left(
1+\frac{m^+(T_a,x)-m^+(T_b,x)}{m^+(T_b,x)}\right)dx\nonumber\\[2mm]
&=&\frac{2}{\pi}\int_0^{+\infty}\frac{1}{x^2}\log\left(1+
\frac{(\Delta-2)x^6
(x^2-(\Delta-2))}{B_1+B_2\frac{m^+(P_{t-4},x)}{m^+(P_{t-3},x)}}\right)dx.
\end{eqnarray}

We first consider the case that $t$ is odd and $t\geq 5$. In the
proof of Theorem \ref{LJtheo1}, we know that the function
$\frac{1}{x^2}\log\left(1+\frac{m^+(T_a,x)-m^+(T_b,x)}
{m^+(T_b,x)}\right)$ is uniformly convergent. Therefore, by Eq.
(\ref{LJequ2011}) and Lemma \ref{LJlem3.5}, we have
\begin{eqnarray*}
&&E(T_a)-E(T_b)\\[2mm]
&>&\frac{2}{\pi}\int_{\sqrt{2}}^{+\infty}\frac{1}{x^2}\log\left(1+
\frac{2x^6(x^2-2)}{B_1+B_2\frac{2}{1+\sqrt{1+4x^2}}}\right)dx+
\frac{2}{\pi}\int_0^{\sqrt{2}}\frac{1}{x^2}\log\left(1+ \frac{2x^6
(x^2-2)}{B_1+B_2\frac{1}{1+x^2}}\right)dx\\[2mm]
&>&\frac{2}{\pi}\cdot 0.02088>0.
\end{eqnarray*}

If $t$ is even, we want to find $t$ and $x$ satisfying that
\begin{eqnarray}\label{LJequ201101}
\frac{m^+(P_{t-4},x)}{m^+(P_{t-3},x)}<\frac{2}{-1+\sqrt{1+4x^2}}.
\end{eqnarray}
It
is equivalent to solve
$$\frac{\lambda_1^{t-3}-\lambda_2^{t-3}}{\lambda_1^{t-2}-
\lambda_2^{t-2}}<-\frac{1}{\lambda_2},$$ which means to solve
$$\left(\frac{\lambda_1}{-\lambda_2}\right)^{t-3}>-2\lambda_2,$$
that is
$$\left(\frac{1+\sqrt{1+4x^2}}{2x}\right)^{2t-6}>\sqrt{1+4x^2}-1.$$
Thus,
$$2t-6>\log_{\frac{1+\sqrt{1+4x^2}}{2x}}(\sqrt{1+4x^2}-1).$$
Since for $x\in (0,+\infty)$, $\frac{1+\sqrt{1+4x^2}}{2x}$ is
decreasing and $\sqrt{1+4x^2}-1$ is increasing, we have that
$\log_{\frac{1+\sqrt{1+4x^2}}{2x}}(\sqrt{1+4x^2}-1)$ is increasing.
Thus, if $x\in [\sqrt{2},5]$, then
$$\log_{\frac{1+\sqrt{1+4x^2}}{2x}}(\sqrt{1+4x^2}-1)\leq
\log_{\frac{1+\sqrt{101}}{10}}(\sqrt{101}-1)<23.$$ Therefore, when
$t\geq 15$, i.e., $2t-6> 23$, we have that Ineq. \eqref{LJequ201101}
holds for $x\in [\sqrt{2},5]$.

Now we calculate the difference of $E(T_a)$ and $E(T_b)$. When $t$
is even and $t\geq 15$, from Eq. \eqref{LJequ2011}, we have
\begin{eqnarray*}
&&E(T_a)-E(T_b)\\[2mm]
&>&\frac{2}{\pi}\int_{5}^{+\infty}\frac{1}{x^2}\log\left(1+
\frac{2x^6(x^2-2)}{B_1+B_2}\right)dx+
\frac{2}{\pi}\int_{\sqrt{2}}^{5}\frac{1}{x^2}\log\left(1+
\frac{2x^6(x^2-2)}{B_1+B_2\frac{2}{-1+\sqrt{1+4x^2}}}\right)dx\\[2mm]
&&+\frac{2}{\pi}\int_0^{\sqrt{2}}\frac{1}{x^2}\log\left(1+
\frac{2x^6
(x^2-2)}{B_1+B_2\frac{2}{1+\sqrt{1+4x^2}}}\right)dx\\[2mm]
&>&\frac{2}{\pi}\cdot 0.003099>0.
\end{eqnarray*}

For $t=3$ and any even $t$ satisfying $4\leq t\leq 14$, by comparing
the energies of the two graphs directly by computer programm, we get
that $E(T_a)<E(T_b)$ for $t=4$, and $E(T_a)>E(T_b)$ for other cases.

The proof is thus complete.\qed

The following theorem gives the result for the cases of
$\Delta=5,6,7$.

\begin{theo}\label{LJtheo4}
For trees with $n$ vertices and two vertices of maximum degree
$\Delta$, let $t=n-4\Delta+4\geq 3$. Then\\
(i) for $\Delta=5$, the maximal energy tree is the graph $T_a$ if
$t$ is odd and $3\leq t\leq89$, and the graph $T_b$ otherwise. \\
(ii) for $\Delta=6$, the maximal energy tree is
the graph $T_a$ if $t=3,5,7$, and the graph $T_b$ otherwise. \\
(iii) for $\Delta=7$, the maximal energy tree is the graph $T_b$ for
any $t\geq 3$.
\end{theo}

\pf In the proof of Theorem \ref{LJtheo1}, we know that the function
$\frac{1}{x^2}\log\left(1+\frac{m^+(T_a,x)-m^+(T_b,x)}
{m^+(T_b,x)}\right)$ is uniformly convergent. We consider the
following cases separately:

(i) $\Delta=5$.

If $t$ is even, we want to find $t$ and $x$ satisfying that
\begin{eqnarray}\label{LJequ201102}
\frac{m^+(P_{t-4},x)}{m^+(P_{t-3},x)}<\frac{2.1}{1+\sqrt{1+4x^2}}.
\end{eqnarray}
It is equivalent to solve
$$\frac{\lambda_1^{t-3}-\lambda_2^{t-3}}{\lambda_1^{t-2}-
\lambda_2^{t-2}}<\frac{2.1}{2\lambda_1},$$ which means to solve
$$\left(\frac{\lambda_1}{-\lambda_2}\right)^{t-3}>\frac{-2.1\lambda_2+2\lambda_1}
{0.1\lambda_1},$$ that is,
$$\left(\frac{1+\sqrt{1+4x^2}}{2x}\right)^{2t-6}>41-\frac{42}{\sqrt{1+4x^2}+1}.$$
Thus,
$$2t-6>\log_{\frac{1+\sqrt{1+4x^2}}{2x}}\left(41-\frac{42}{\sqrt{1+4x^2}+1}
\right).$$ Since for $x\in (0,+\infty)$,
$\frac{1+\sqrt{1+4x^2}}{2x}$ is decreasing and
$-\frac{42}{\sqrt{1+4x^2}+1}$ is increasing, we have that
$\log_{\frac{1+\sqrt{1+4x^2}}{2x}}\left(41-\frac{42}{\sqrt{1+4x^2}+1}
\right)$ is increasing. Thus, if $x\in (0,\sqrt{3}]$,
$$\log_{\frac{1+\sqrt{1+4x^2}}{2x}}\left(41-\frac{42}{\sqrt{1+4x^2}+1}
\right)\leq
\log_{\frac{1+\sqrt{13}}{2\sqrt{3}}}\left(41-\frac{42}{1+\sqrt{13}}\right)<13.$$
Therefore, when $t\geq 10$, i.e., $2t-6>13$, we have that Ineq.
\eqref{LJequ201102} holds for $x\in (0,\sqrt{3}]$. Thus, if $t$ is
even and $t\geq 10$, from Eq. \eqref{LJequ2011} and Lemma
\ref{LJlem3.5}, we have
\begin{eqnarray*}
E(T_a)-E(T_b)
&<&\frac{2}{\pi}\int_{\sqrt{3}}^{+\infty}\frac{1}{x^2}\log\left(1+
\frac{3x^6(x^2-3)}{B_1+B_2\frac{2}{1+\sqrt{1+4x^2}}}\right)dx\\[2mm]
&&+ \frac{2}{\pi}\int_{0}^{\sqrt{3}}\frac{1}{x^2}\log\left(1+
\frac{3x^6(x^2-3)}{B_1+B_2\frac{2.1}{1+\sqrt{1+4x^2}}}\right)dx\\[2mm]
&<&\frac{2}{\pi}\cdot(-4.43\times10^{-4})<0.
\end{eqnarray*}

If  $t$ is odd, we want to find $t$ and $x$ satisfying that
\begin{eqnarray}\label{LJequ201103}
\frac{m^+(P_{t-4},x)}{m^+(P_{t-3},x)}>\frac{1.99}{1+\sqrt{1+4x^2}},
\end{eqnarray} that is
$$2t-6>\log_{\frac{1+\sqrt{1+4x^2}}{2x}}\left(399-\frac{398}{\sqrt{1+4x^2}+1}
\right).$$ Since for $x\in(0,+\infty)$,
$\log_{\frac{1+\sqrt{1+4x^2}}{2x}}\left(399-\frac{398}{\sqrt{1+4x^2}+1}
\right)$ is increasing, we have that if $x\in [\sqrt{3},390]$,
$$\log_{\frac{1+\sqrt{1+4x^2}}{2x}}\left(399-\frac{398}{\sqrt{1+4x^2}+1}
\right)<4671.$$ Therefore, for $t\geq 2339$, i.e., $2t-6\geq 4671$,
we have that Ineq. \eqref{LJequ201103} holds for $x\in
[\sqrt{3},390]$. Thus, if $t$ is odd and $t\geq 2339$, from Eq.
\eqref{LJequ2011} and Lemma \ref{LJlem3.5}, we have
\begin{eqnarray*}
&&E(T_a)-E(T_b)\\[2mm]
&<&\frac{2}{\pi}\int_{390}^{+\infty}\frac{1}{x^2}\log\left(1+
\frac{3x^6(x^2-3)}{B_1+B_2\frac{1}{1+x^2}}\right)dx
+\frac{2}{\pi}\int_{\sqrt{3}}^{390}\frac{1}{x^2}\log\left(1+
\frac{3x^6(x^2-3)}{B_1+B_2\frac{1.99}{1+\sqrt{1+4x^2}}}\right)dx\\[2mm]
 &&+\frac{2}{\pi}\int_{0}^{\sqrt{3}}\frac{1}{x^2}\log\left(1+
\frac{3x^6(x^2-3)}{B_1+B_2\frac{2}{1+\sqrt{1+4x^2}}}\right)dx\\[2mm]
&<&\frac{2}{\pi}\cdot(-6.66\times10^{-6})<0.
\end{eqnarray*}

For any even $t$ satisfying that $4\leq t\leq 8$ and any odd $t$
satisfying that $3\leq t\leq 2337$, by comparing the energies of the
two graphs directly by matlab programm, we get that $E(T_a)>E(T_b)$
for any odd $t$ satisfying $3\leq t\leq89$, and $E(T_a)<E(T_b)$ for
the other cases.

(ii) $\Delta=6$.

If $t$ is even and $t\geq 4$, from Eq. \eqref{LJequ2011} and Lemma
\ref{LJlem3.5}, we have
\begin{eqnarray*}
E(T_a)-E(T_b)
&<&\frac{2}{\pi}\int_{2}^{+\infty}\frac{1}{x^2}\log\left(1+
\frac{4x^6(x^2-4)}{B_1+B_2\frac{2}{1+\sqrt{1+4x^2}}}\right)dx\\[2mm]
&&+ \frac{2}{\pi}\int_{0}^{2}\frac{1}{x^2}\log\left(1+
\frac{4x^6(x^2-4)}{B_1+B_2}\right)dx\\[2mm]
&<&\frac{2}{\pi}\cdot(-0.02027)<0.
\end{eqnarray*}

If $t$ is odd, similar to the proof in (i), we can show that when
$t\geq 27$ and $x\in [2,22]$, the following inequality holds:
$$\frac{m^+(P_{t-4},x)}{m^+(P_{t-3},x)}>\frac{1}{1+\sqrt{1+4x^2}}.$$
Hence, if $t$ is odd and $t\geq 27$, we have
\begin{eqnarray*}
&&E(T_a)-E(T_b)\\[2mm]
&<&\frac{2}{\pi}\int_{22}^{+\infty}\frac{1}{x^2}\log\left(1+
\frac{4x^6(x^2-4)}{B_1+B_2\frac{1}{1+x^2}}\right)dx
+\frac{2}{\pi}\int_{2}^{22}\frac{1}{x^2}\log\left(1+
\frac{4x^6(x^2-4)}{B_1+B_2\frac{1}{1+\sqrt{1+4x^2}}}\right)dx\\[2mm]
 &&+\frac{2}{\pi}\int_{0}^{2}\frac{1}{x^2}\log\left(1+
\frac{4x^6(x^2-4)}{B_1+B_2\frac{2}{1+\sqrt{1+4x^2}}}\right)dx\\[2mm]
&<&\frac{2}{\pi}\cdot(-2.56\times10^{-4})<0.
\end{eqnarray*}

For any odd $t$ satisfying that $3\leq t\leq 25$, by comparing the
energies of the two graphs directly by matlab programm, we get that
$E(T_a)>E(T_b)$ for $t=3,5,7$, and $E(T_a)<E(T_b)$ for the other
cases.

(iii) $\Delta=7$.

If $t$ is even and $t\geq 4$, by the same method as used in (ii), we
get that $E(T_a)-E(T_b)<\frac{2}{\pi}\cdot(-0.04445)<0.$

If $t$ is odd and $t\geq 5$, we have that
\begin{eqnarray*}
E(T_a)-E(T_b)&<&\frac{2}{\pi}\int_{\sqrt{5}}^{+\infty}\frac{1}{x^2}\log\left(1+
\frac{5x^6(x^2-5)}{B_1+B_2\frac{1}{1+x^2}}\right)dx\\[2mm]
&&+ \frac{2}{\pi}\int_{0}^{\sqrt{5}}\frac{1}{x^2}\log\left(1+
\frac{5x^6(x^2-5)}{B_1+B_2\frac{2}{1+\sqrt{1+4x^2}}}\right)dx\\[2mm]
&<&\frac{2}{\pi}\cdot(-0.01031)<0.
\end{eqnarray*}

For $t=3$, we can compare the energies of the two graphs directly by
matlab programm and get that $E(T_a)<E(T_b)$.

The proof is now complete.\qed

\end{document}